\documentclass[review]{elsarticle}
\usepackage{lineno}
\usepackage{bbold}
\usepackage{verbatim}
\usepackage{url}
\usepackage{graphicx}
\usepackage{float}
\usepackage{mathrsfs}
\usepackage{epstopdf}
\usepackage[all,color]{xy}
\usepackage{cancel}
\usepackage{MnSymbol}
\usepackage{tikz}
\usepackage[margin=1in]{geometry}

\usetikzlibrary{arrows}
\modulolinenumbers[5]
\allowdisplaybreaks
\newtheorem{theorem}{Theorem}

\newtheorem{corollary}[theorem]{Corollary}
\newtheorem{example}[theorem]{Example}
\newtheorem{lemma}[theorem]{Lemma}

\newproof{proof}{Proof}
\numberwithin{equation}{subsection}
\numberwithin{theorem}{subsection}

\newcommand{\ignore}[1]{}

\DeclareMathOperator{\sgn}{sgn}

\begin{document}

\begin{frontmatter}

\title{Generalizing Kirchhoff laws for Signed Graphs}

\author[add2]{Lucas J. Rusnak\corref{mycorrespondingauthor}}\ead{Lucas.Rusnak@txstate.edu}
\author[add2]{Josephine Reynes} %jar369@txstate.edu
\author[add1]{Skyler J. Johnson} %dawsonjohnson314@gmail.com
\author[add1]{Peter Ye} %peterye03@gmail.com

\address[add2]{Department of Mathematics, Texas State University, San Marcos, TX 78666, USA}

\address[add1]{Mathworks, Texas State University, San Marcos, TX 78666, USA}

\cortext[mycorrespondingauthor]{Corresponding author}

\begin{abstract}

Kirchhoff-type Laws for signed graphs are characterized by generalizing transpedances through the incidence-oriented structure of bidirected graphs. The classical $2$-arborescence interpretation of Tutte is shown to be equivalent to single-element Boolean classes of reduced incidence-based cycle covers, called contributors. A generalized contributor-transpedance is introduced using entire Boolean classes that naturally cancel in a graph; classical conservation is proven to be property of the trivial Boolean classes. The contributor-transpedances on signed graphs are shown to produce non-conservative Kirchhoff-type Laws, where every contributor possesses the unique source-sink path property.  Finally, the maximum value of a contributor-transpedance is calculated through the signless Laplacian.

\end{abstract}


\end{frontmatter}

%%%%%%%%

\section{Introduction and Background}

We introduce and characterize Kirchhoff-type Laws for signed graphs by generalizing transpedances from \cite{BSST}. This is accomplished by using the incidence-theoretic approach introduced in \cite{AH1} to study hypergraphic Laplacians, and the incidence-path mapping families, called \emph{contributors}, from \cite{OHSachs} that generalize cycle covers to classify various hypergraphic characteristic polynomials similar to Sachs' Theorem \cite{Sim1,SGBook}. It was shown in \cite{OHMTT} that if all edges are size $2$ these generalized cycle covers form Boolean lattices that generalize the Matrix-tree theorem. These Boolean families are naturally cancellative when $G$ is a graph with the trivial single-element classes corresponding to spanning trees --- providing ``conservation'' for the graphic Kirchhoff Laws. 

Transpedances were introduced in \cite{BSST} as a way to study the packing and cutting problem of dissecting a rectangle into squares by translating the question into a networking potential problem. A graph is associated to each dissection and its natural flow capacity is determined to be the tree-number via the Matrix-tree Theorem. Moreover, the size of the admissible squares are ordered second-cofactors of the Laplacian, and a combinatorial interpretation of Kirchhoff's Laws via ``spanning tree flows'' is obtained for any source-sink pair where edges are labeled by signed $2$-arborescences. A brief introduction to transpedances appears in Subsection \ref{ssec:TutteIntro}. The original investigation into transpedances is credited with leading to Tutte's investigation into Graph-polynomials \cite{tutte2004graph}. Non-conservative Kirchhoff-type Laws for directed graphs appear in \cite{Tutte}, the algebraic development of potential theory appears in \cite{biggs_1997}, and a formulation using ported matroids appear in \cite{chaiken2006ported} that analyze signed contributions of spanning forests --- in this paper, we identify non-forest contributors in Boolean equivalence classes that produce non-conservative generalizations of Kirchhoff's Laws.

A signed graph is a generalization of a graph where each edge receives a sign $+1$ or $-1$ to examine social balance \cite{Har0}. Incidence-orientations of signed graphs \cite{OSG}, or bidirected graphs, first appeared in integer programming \cite{MR0267898}. This incidence-theoretic approach has led to the incidence-oriented hypergraphic characterization of the Laplacian \cite{OHSachs,AH1,OHMTT} that generalizes the signed graphic All Minors Matrix-tree theorem \cite{Seth1} as well as the signed graphic Sachs' Theorem \cite{Sim1}. An introduction to the necessary incidence-theoretic concepts appear in Subsection \ref{ssec:OHIntro}.

A simple interpretation for contributor families for any oriented hypergraph is introduced in Section \ref{sec:ContClass} that specializes to the Boolean lattice equivalence classes in \cite{OHMTT}. Ordered second principal order ideals correspond to the ordered second cofactors that determine the coefficients of the total-minor polynomial in \cite{IH2}. The single element Boolean classes are shown to be in bijective correspondence with the $2$-arborescences from \cite{BSST} by using the Linking Lemma. The notion of transpedances is extended in Section \ref{sec:Transp} to \emph{$D$-contributor-transpedances}, which include all the Boolean classes, not just the single element ones. Kirchhoff's Degeneracy and Energy Reversal conditions are shown to immediately hold for them. The $D$-contributor-transpedance value is then calculated for an arbitrary edge and the Boolean classes are shown to vanish if they contain a positive circle. Thus, if $G$ is a graph, then only the trivial contributors that correspond to $2$-arborescences remain.

Section \ref{sec:Kirch} proves that all $D$-contributor-transpedances possess a unique source-sink path property, and the trivial classes used to label the edges sort spanning trees along their source-sink path. Kirchhoff's Cycle and Vertex Conservation Laws are shown to be a property of the trivial Boolean classes, and conservation on non-cancellative Boolean classes (negative classes) cannot be guaranteed. The maximal contributor-transpedance problem is solved in Section \ref{sec:MaxPerm} through the signless Laplacian and the permanent. This count holds for any oriented hypergraph, and a simple permanent version of Kirchhoff's Laws via contributors is stated. Unfortunately, the techniques for a complete general hypergraphic transpedance version are limited by: (1) the partial parallel edges that readily appear, and (2) the lack of Boolean nature of general classes. We hope the general class partial order introduced in Section \ref{sec:ContClass} may serve to remedy this. Additionally, the connection to Tutte-functions and ported matroids require further study to relate to the work in \cite{CHAIKEN198996, chaiken2006ported}.

\subsection{Transpedances and Tutte's Results}
\label{ssec:TutteIntro}

A \emph{$2$-arborescence} of $G$ is a pair of disjoint rooted trees whose union spans $G$. Kirchhoff's Laws with unit resistance has been shown to be equivalent to $2$-arborescence counts whose values are commensurable with the tree number of the graph \cite{BSST}; non-unit resistance is simply a weighted version of this combinatorial result, while directed graphs produce a non-conservative version of Kirchhoff's Laws \cite{Tutte}. Let $u_1, u_2, w_1, w_2 \in V(G)$, and define $\langle u_1 w_1, u_2 w_2 \rangle$ be the number of $2$-arborescences with one component rooted at $u_1$ and containing $w_1$,  and the other component rooted at $u_2$ and containing $w_2$.

\begin{figure}[H]
    \centering
    \includegraphics{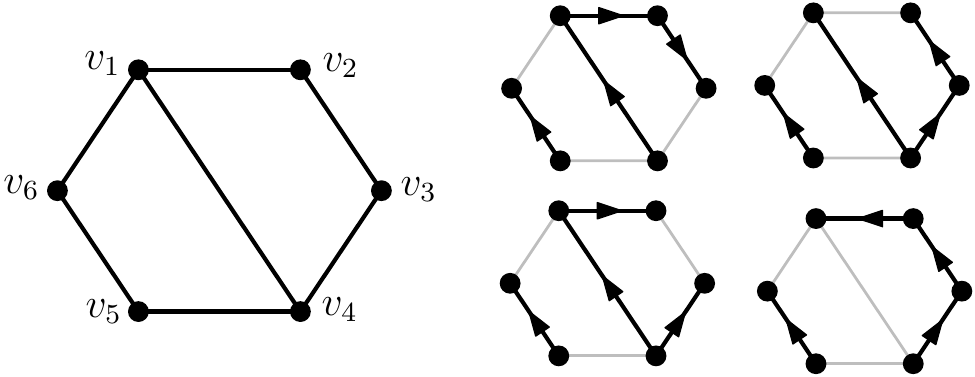}
    \caption{All $2$-arborescences of the graph of the form $\langle v_5 v_6, v_4 v_1 \rangle$.}
    \label{fig:2-arbors}
\end{figure}

Given graph $G$ with source $u_1$ and sink $u_2$, the \emph{$w_1 w_2$-transpedance} of $G$ is 
\begin{align*}
    [u_1 u_2, w_1 w_2] = \langle u_1 w_1, u_2 w_2 \rangle - \langle u_1 w_2, u_2 w_1 \rangle.
\end{align*}
It was shown in \cite{BSST,Tutte} that the value $[u_1 u_2, w_1 w_2]$ is also the (ordered) second cofactor of the Laplacian of $G$. Let $\mathbf{L}_{G}$ be the Laplacian of $G$, let $\mathbf{L}_{(G;u_1,w_1)}$ be the $u_1 w_1$-minor of $L$, let $\mathbf{L}_{(G;u_1 u_2,w_1 w_2)}$ be the $u_2 w_2$-minor of $\mathbf{L}_{(G;u_1,w_1)}$, and define $\mathbf{L}_{(G;\mathbf{u},\mathbf{w})}$ iteratively for vertex vectors $\mathbf{u},\mathbf{w}$. Specifically, $[u_1 u_2, w_1 w_2]$ is the value of the $u_2 w_2$-cofactor in the $u_1 w_1$-minor using the positional sign of $u_1 w_1$ in $\mathbf{L}_{G}$ and the positional sign of $u_2 w_2$ in $\mathbf{L}_{(G;u_1,w_1)}$.

\begin{example}
Since, in Figure \ref{fig:2-arbors}, there are no $2$-arborescences of the form $\langle v_5 v_1, v_4 v_6 \rangle$, the transpedance value $[v_5 v_4, v_6 v_1] = 4$ is assigned to the edge between $v_5$ and $v_6$. Note that transpedances are directional, so $[v_5 v_4, v_6 v_1]$ can be regarded as the potential drop from $v_6$ to $v_1$ with source $v_5$ and sink $v_4$. Thus, $[v_5 v_4, v_1 v_6]$ would be $-4$, but would arise from a different the set of $2$-arborescences.
\end{example}

Edge labeling by transpedances produces a combinatorial Kirchhoff's Laws that are summarized as follows:

\begin{theorem}[\cite{BSST,Tutte}]
\label{t:TutteTransp}
Let $G$ be a graph with tree number $\tau(G)$, the following hold:
\begin{enumerate}
    \item (Degeneracy) $[u_1 u_1, w_1 w_2] = [u_1 u_2, w_1 w_1] = 0$,
    \item (Energy Reversal) $[u_1 u_2, w_1 w_2] = -[u_1 u_2, w_2 w_1]= -[u_2 u_1, w_1 w_2]$,
    \item (Cycle Conservation) $[u_1 u_2, w_1 w_2] + [u_1 u_2, w_2 w_3] + [u_1 u_2, w_3 w_1] = 0$,
    \item (Vertex Conservation) $\displaystyle
    \dsum_{y:y \sim w_1} l_{vy}[u_1 u_2, w_1 y] = \tau(G) \delta_{u_2 w_1} - \tau(G) \delta_{u_1 w_1}$,
\end{enumerate}
where $\delta_{u w} = 1$ if $u=w$, and is $0$ otherwise.
\end{theorem}
Part (1) establishes that degenerate transpedances have a value of $0$, part (2) is reversal of flow, (3) implies both path concatenation and cycle-conservation, and (4) is vertex-conservation, with the exception of the source and sink where the edges have natural flow of $\tau(G)$ out of the source, and into the sink. 

\begin{example}
The transpedance labeling of the graph in Figure \ref{fig:2-arbors} with source $v_5$ and sink $v_4$ appears in Figure \ref{fig:KirchA}. The four $2$-arborescences in Figure \ref{fig:2-arbors} are assigned to the directed adjacency between $v_6$ and $v_1$.
\begin{figure}[H]
    \centering
    \includegraphics{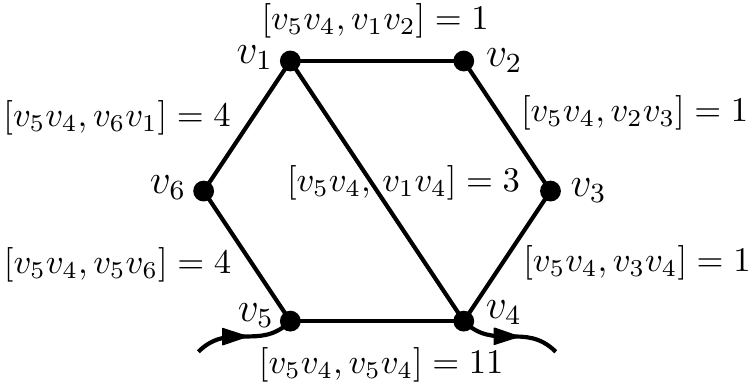}
    \caption{A transpedance labeling of $G$ with source $v_5$ and sink $v_4$}
    \label{fig:KirchA}
\end{figure}
There are $15$ spanning trees, represented as the net inflow and outflow from $v_5$ and $v_4$, respectively. It is easy to check that the directed cycle sums relative to source $v_5$ and sink $v_4$ are zero. Also, the in/out vertex sums are also zero --- with the exception of the source and sink, whose values are the tree number $15$.
\end{example}

\subsection{Incidence Orientations and Signed Graphs}
\label{ssec:OHIntro}

An incidence hypergraph is a quintuple $G=(V, E, I, \varsigma, \omega)$ consisting of a set of vertices $V$, a set of edges $E$, a set of incidences $I$, and two incidence maps $\varsigma: I \to V$ and $\omega: I \to E$. An \emph{orientation of an incidence hypergraph G} is a signing function $\sigma:I\rightarrow\{+1,-1\}$, which produces a $V \times E$ integer incidence matrix $\mathbf{H}_{G}$. The \emph{Laplacian matrix of $G$} is defined as $\mathbf{L}_{G}:=\mathbf{H}_{G} \mathbf{H}_{G}^{T}=\mathbf{D}_{G}-\mathbf{A}_{G}$, where the degree matrix is the number of incidences at a vertex, and the adjacency matrix has entries determined by the sign $-\sigma(i)\sigma(j)$, where $i$ and $j$ are the incidences of an adjacency \cite{AH1}. A \emph{bidirected graph} is an oriented hypergraph in which every edge is a $2$-edge. Bidirected graphs first appeared in integer programming \cite{MR0267898}, and later were shown to be orientations of signed graphs \cite{OSG}. The edge labeling of a bidirected graph by the adjacency sign  $-\sigma(i)\sigma(j)$ is called a \emph{signed graph}, and a graph can be regarded as a signed graph with all edges positive.
\begin{figure}[H]
    \begin{align*}
        \vcenter{\hbox{
        \includegraphics{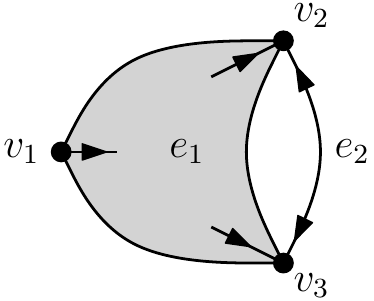}}}
    & &
        \mathbf{L}_{G} = \left[\begin{array}{ccc}
        1 & -1 & -1  \\ 
        -1 & 2 & 2  \\
        -1 & 2 & 2 
        \end{array}\right]
    \end{align*}
    \caption{An oriented hypergraph $G$ and its Laplacian.}
    \label{fig:OHLap}
\end{figure}
A \emph{contributor of $G$} is an incidence preserving map from a disjoint union of $\overrightarrow{P}_{1}$'s with tail $t$ and head $h$ into $G$ defined by $c:\dcoprod \limits_{v\in V}\overrightarrow{P}_{1}\rightarrow G$ such that $c(t_{v})=v$ and $\{c(h_{v})\mid v\in V\}=V$. Due to the nature of the incidence-maps it is possible for a path to fold back on itself creating a \emph{backstep} of the form $v,i,e,i,v$ --- these are the entries in the hypergraphic degree matrix. A contributor can be regarded as a permutation clone that is a generalized cycle covers similar to Sachs' Theorem to determine characteristic polynomial coefficients \cite{SGBook,Sim1,OHSachs}; contributors naturally form Boolean lattices when $G$ is a bidirected graph \cite{OHMTT}. The set of contributors of an oriented hypergraph is denoted ${\mathfrak{C}}(G)$. Throughout, let $U,W \subseteq V$ with $\lvert U \rvert = \lvert W \rvert$, while a total ordering of each set will be denoted by $\mathbf{u}$ and $\mathbf{w}$, respectively. Let ${\mathfrak{C}}(G;\mathbf{u},\mathbf{w})$ be the set of \emph{restricted} contributors in $G$ where $c(u_i)=w_i$, and two elements of ${\mathfrak{C}}(G;\mathbf{u},\mathbf{w})$ are said to be \emph{$[\mathbf{u},\mathbf{w}]$-equivalent}. Let $\widehat{\mathfrak{C}}(G;\mathbf{u},\mathbf{w})$ be the set obtained by removing the $\mathbf{u} \rightarrow \mathbf{w}$ mappings from ${\mathfrak{C}}(G;\mathbf{u},\mathbf{w})$; the elements of $\widehat{\mathfrak{C}}(G;\mathbf{u},\mathbf{w})$ are called the \emph{reduced} $[\mathbf{u},\mathbf{w}]$-equivalent contributors. To avoid confusion between an algebraic cycle and a graph component that forms a closed walk we refer to the graph images as \emph{circles}, and backsteps will be considered separate from circles as they do not complete an adjacency.

\begin{example}
Figure \ref{fig:OHContrib} shows four contributors of the hypergraph $G$ from Figure \ref{fig:OHLap}. The tail of each path is labeled with a different shape and mapped to its corresponding vertex in $G$; the heads are then mapped to again cover the vertices. The bottom two contributors consist of all backsteps and are clones of the identity permutation (but are distinct contributors). 
\begin{figure}[H]
    \centering
    \includegraphics{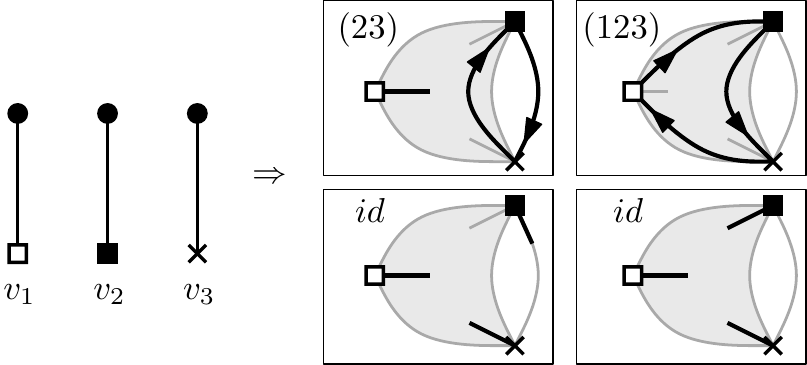}
    \caption{Contributor examples with associated permutations.}
    \label{fig:OHContrib}
\end{figure}
The contributors above each identity-clone are  $[v_2,v_3]$-equivalent as there is a path mapping to the $v_2 v_3$-adjacency in each contributor.
\end{example}

Transpedances are second cofactors and these arise naturally as the coefficients of the degree-$2$ monomials of the total minor polynomial for integer matrix Laplacians \cite{IH2}, where the coefficients are determined by sums of reduced contributors. Let $\chi ^{P}(\mathbf{L}_{G},\mathbf{x})$ and $\chi ^{D}(\mathbf{L}_{G},\mathbf{x})$ be the total minor polynomial as determined by the permanent and determinant of $\mathbf{X} - \mathbf{L}_{G}$, respectively, where the $ij$-entry of $\mathbf{X}$ is $x_{ij}$. The total minor polynomial is calculated as follows:

\begin{theorem}[\cite{IH2}, Theorem 3.1.2]
\label{t:Poly}
Let $G$ be an oriented hypergraph with Laplacian matrix $\mathbf{L}_{G}$, then \begin{enumerate}
\item $\chi ^{P}(\mathbf{L}_{G},\mathbf{x})=\dsum\limits_{ [\mathbf{u},\mathbf{w}]}
\left(
\dsum\limits_{\substack{c \in \widehat{\mathfrak{C}}(L^{0}(G);\mathbf{u},\mathbf{w}) \\ \sgn(c) \neq 0}}
(-1)^{nc(c)+bs(c)}
\right)
\dprod\limits_{i}x_{u_{i},w_{i}}$,

\item $\chi ^{D}(\mathbf{L}_{G},\mathbf{x})=\dsum\limits_{ [\mathbf{u},\mathbf{w}]}
\left(
\dsum\limits_{\substack{c \in \widehat{\mathfrak{C}}(L^{0}(G);\mathbf{u},\mathbf{w}) \\ \sgn(c) \neq 0}}
(-1)^{ec(\check{c})+nc(c)+bs(c)}
\right)
\dprod\limits_{i}x_{u_{i},w_{i}}$.
\end{enumerate}
where $ec(\check{c})$ represents the number of even-cycles in the unreduced contributor of $c$, $bs(c)$ represents the number of backsteps, and $nc(c)$ represents the number of negative components.
\end{theorem} 
The hypergraph $L^{0}(G)$ is the \emph{zero-loading of $G$} and extends the hypergraph to a uniform hypergraph and assigns a weight of $0$ to all new incidences. Thus, a reduced contributor exists in $G$ if and only if it is non-zero. As discussed in \cite{IH2}, the $\mathbf{u} \rightarrow \mathbf{w}$ maps need not exist in $G$ as they are removed, but the maps must be allowed to exist a priori their removal, which is remedied by the zero-loading $L^{0}(G)$. To simplify notation let $\widehat{\mathfrak{C}}_{\neq0}(L^{0}(G);\mathbf{u},\mathbf{w})$ be the set of non-zero reduced contributors in $L^{0}(G)$; that is, the reduced contributors that reside in $G$.

\begin{example}
Consider the value $[v_5 v_4, v_1 v_2] = 1$ along the top edge in Figure \ref{fig:KirchA}. To find this value using the total minor polynomial we first find all contributors where $v_5 \mapsto v_1$ and $v_4 \mapsto v_2$, then remove these two maps --- these maps are allowed to exist in the zero-loading $L^{0}(G)$ and are subsequently removed, the remaining objects need to exist in $G$ to avoid mapping to $0$. There is only one such reduced contributor that lies in $G$, shown in Figure \ref{fig:TotPolyEx}.

\begin{figure}[H]
    \centering
    \includegraphics{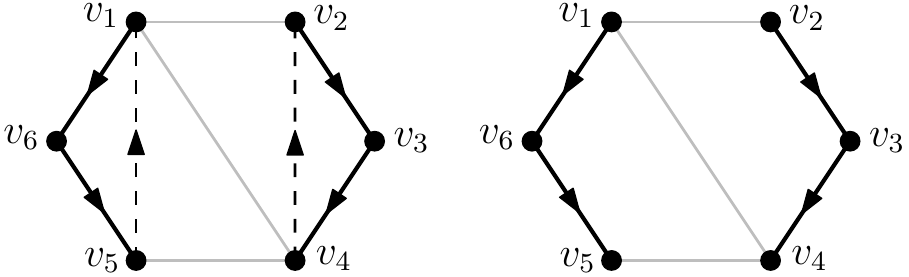}
    \caption{Reduced contributors find coefficients of the total minor polynomial as generalized cycle covers.}
    \label{fig:TotPolyEx}
\end{figure}
Using the determinant signing function in Theorem \ref{t:Poly}, and assuming every edge is positive (as in a graph), we have $ec(\check{c}) = 0$ since there are $0$ even circles in the non-reduced contributor, while $nc(c)=0$ and $bs(c) = 0$. Thus, the sign of the contributor is $(-1)^{0+0+0}=1$. This is the value of the coefficient of $x_{v_5 v_1}x_{v_4 v_2}$ as well as $[v_5 v_4, v_1 v_2]$, as depicted in Figure \ref{fig:KirchA}.
\end{example}

We show that the contributor mappings produce a natural adjacency labeling for oriented hypergraphs and a non-conservative generalization of Kirchhoff's Laws for signed graphs. Moreover, this process is not limited to the determinant. If every adjacency is negative, then the permanent counts the total number of contributors for an edge; thus, providing a maximum value for potential on each edge.

\section{Contributor Classes and Arborescences}
\label{sec:ContClass}

\subsection{Tail Equivalence}
\label{ssec:TailEquiv}

\begin{figure}[H]
    \centering
    \includegraphics{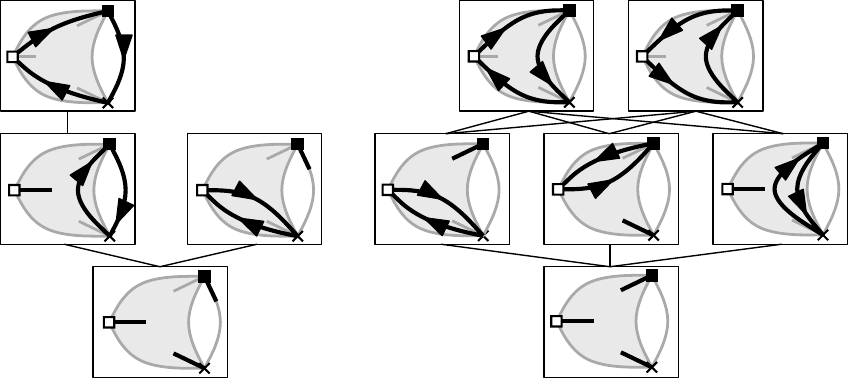}
    \caption{Tail-equivalence classes from Figure \ref{fig:OHContrib}.}
    \label{fig:TailEq}
\end{figure}

Two contributors are said to be \emph{tail-equivalent} if the image of their tail-incidences agree. Each identity-contributor in Figure \ref{fig:OHContrib} is tail-equivalent to the contributor above it, as they both enter the same edge, but complete to different permutations. Clearly, there is exactly one identity-contributor in each tail-equivalency class. The elements of a tail-equivalence class are partially ordered by $c \leq c'$ if (1) the set of circles of $c$ is contained in the set of circles of $c'$, or (2) the set of incidences are equal and $c$ has more connected components than $c'$. Thus, the identity-contributor, having the most components and an empty set of circles, is the least element of each poset, while the number of contributors on a single $k$-edge follow the Stirling numbers of the first kind. Two examples appear in Figure \ref{fig:TailEq}.

The concept of tail-equivalence is a generalization of circle activation classes of bidirected graphs in \cite{OHMTT}, where extending a backstep into its unique directed adjacency was called \emph{unpacking}, and folding a directed adjacency back into a backstep was called \emph{packing}. These operations are well-defined and inverses in a bidirected graph, while in an oriented hypergraph only packing is well-defined on larger edges. Contributors that were packing/unpacking equivalent were grouped into \emph{activation classes} and ordered as new circles appear. 

We adopt the convention of referring to a tail-equivalence class as an activation class when $G$ is a bidirected graph. Let $\mathcal{A}(G)$ denote a tail-equivalence class of $G$. As with restricted and reduced contributors we let $\mathcal{A}(\mathbf{u};\mathbf{w};G)$ be the elements of tail-equivalency class $\mathcal{A}(G)$ where $u_i \mapsto w_i$, and $\hat{\mathcal{A}}(\mathbf{u};\mathbf{w};G)$ be the elements of $\mathcal{A}(\mathbf{u};\mathbf{w};G)$ with $u_{i} \mapsto w_{i}$ removed for each $i$. From \cite{OHMTT}, the activation classes and their restricted subclasses (order ideals) of a bidirected graph are Boolean. 

\begin{lemma}[\cite{OHMTT}, Lemma 3.6]
\label{l:booleanLattice}
For a bidirected graph $G$, all activation classes of $G$ are Boolean lattices.
\end{lemma}
It was also shown in \cite{IH2} that the reduced contributors in single element activation classes $\hat{\mathcal{A}}_{\neq 0}(\mathbf{u};\mathbf{w};L^0(G))$ are unpacking equivalent to $k$-arborescences.

\begin{theorem}[\cite{IH2}, Theorem 3.2.4]
\label{t:karbor}
In a bidirected graph $G$ the set of all elements in single-element $\hat{\mathcal{A}}_{\neq 0}(\mathbf{u};\mathbf{w};L^{0}(G))$ is unpacking equivalent to $k$-arborescences. Moreover, the $i^{th}$ component in the arborescence has sink $u_i$, and the vertices of each component are determined by the linking induced by $c^{-1}$ between all $u_i \in U \cap \overline{W} \rightarrow \overline{U}$ or unpack into a vertex of a linking component.
\end{theorem}

\begin{example}
Consider the graph from Figure \ref{fig:2-arbors} as an incidence-graph. Each identity-contributor has no circles and the backsteps may be unpacked to produce new cycles. Since every edge contains a unique adjacency, the contributors are ordered by their circle sets. Moreover, the subclass of where $v_i \mapsto v_j$ is an order ideal. Three activation classes appear in Figure \ref{fig:ActivationClass} along with their $v_5 \mapsto v_4$ subclasses highlighted. The top contributor in the rightmost figure in Figure \ref{fig:ActivationClass} is a trivial $v_5 \mapsto v_4$ subclass. Additionally, the removal of the $v_5 \mapsto v_4$ map leaves a rooted  spanning tree ($1$-arborescence).

\begin{figure}[H]
\begin{center}
\includegraphics[]{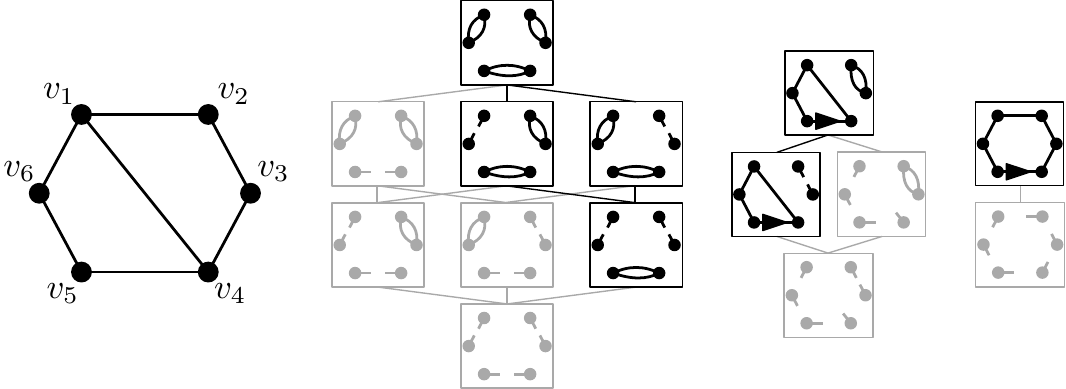}
\end{center}
    \caption{Three Boolean activation classes for the given graph and their $v_5 \mapsto v_4$ activation subclass (darker).}
    \label{fig:ActivationClass}
\end{figure}

To see how a $2$-arborescence is formed consider the middle activation class in Figure \ref{fig:ActivationClass} where $v_5 \mapsto v_4$. Remove the $v_5 \mapsto v_4$ map and then take the second order ideal induced by $v_2 \mapsto v_3$ --- this gives the middle figure in Figure \ref{fig:ReducedActivationClass}. The removal of the $v_2 \mapsto v_3$ mapping (and unpacking any backsteps) yields the $2$-arborescence on the right of Figure \ref{fig:ReducedActivationClass}.

\begin{figure}[H]
\begin{center}
\includegraphics[]{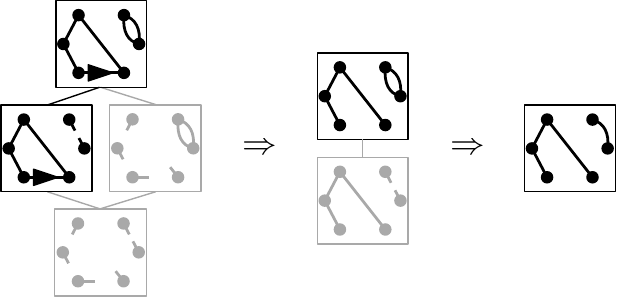}
\end{center}
    \caption{A trivial $[v_5 v_2,v_4 v_3]$-reduced activation class unpacks into a $2$-arborescence.}
    \label{fig:ReducedActivationClass}
\end{figure}

\end{example}

We show that Tutte's transpedances are actually statements about trivial, single-element, activation classes. Since single-element trivial activation classes will appear repeatedly, let $\widehat{\mathfrak{C}}^1_{\neq0}(L^0(G);u_1 u_2, w_1 w_2)$ be the non-zero elements of $\widehat{\mathfrak{C}}(L^0(G);u_1 u_2, w_1 w_2)$ in trivial activation classes.

\subsection{Contributor Arborescences}
\label{ssec:Arbors}

The $2$-arborescences that arise from trivial activation classes need not be the same as Tutte's. A $2$-arborescence for the transpedance calculation $[u_1 u_2, w_1 w_2]$ will be called a \emph{Tutte-$2$-arborescence}, while a $2$-arborescence described as an element of $\widehat{\mathfrak{C}}^1_{\neq0}(L^0(G);u_1 u_2, w_1 w_2)$ will be called a \emph{contributor-$2$-arborescence}. 

Let $F$ be a Tutte-$2$-arborescence in the calculation of $[u_1 u_2, w_1 w_2]$; the sign of $F$ (relative to $[u_1 u_2, w_1 w_2]$), denoted $sgn_T(F)$, is $+1$ if it contributes to the value of $\langle u_1 w_1, u_2 w_2 \rangle$ and $-1$ if it contributes to the value of $\langle u_1 w_2, u_2 w_1 \rangle$. Tutte and contributor-$2$-arborescences are related via the Linking Lemma and the number of cycles that are formed.

\begin{lemma}
\label{l:ArborBijection}
There is a bijection between Tutte-$2$-arborescences of the form $[u_1 u_2, w_1 w_2]$ and contributor-$2$-arborescences from $\widehat{\mathfrak{C}}^1_{\neq0}(L^0(G);u_1 u_2, w_1 w_2)$. 
\end{lemma}

\begin{proof}
Let $u_1$ and $u_2$ be the source and sink, respectively, and let $w_1$ and $w_2$ be two vertices. 

\textbf{Part I:} Let $F$ be a Tutte-$2$-arborescence for $[u_1 u_2, w_1 w_2]$. There are two cases based on $sgn_T(F)$.

\textit{Case 1 $(sgn_T(F) = +1)$:}  If $sgn_T(F) = +1$, then $u_1$ and $w_1$ are in one component, and $u_2$ and $w_2$ are in the other. Reverse the path from $u_1$ and $w_1$ and $u_2$ and $w_2$ within each component. Introduce edges directed $u_1 \mapsto w_1$ and $u_2 \mapsto w_2$ to complete two disjoint cycles. Note that these edges need not exist in $G$ as they exist in the $0$-loading and will be removed in the reduced contributor. Next, pack all adjacencies away from each cycle into backsteps and remove the $u_1 \mapsto w_1$ and $u_2 \mapsto w_2$ adjacencies. Since there are no more circles, the resulting object is in  $\widehat{\mathfrak{C}}^1_{\neq0}(L^0(G);u_1 u_2, w_1 w_2)$.

\textit{Case 2 $(sgn_T(F) = -1)$:} If $sgn_T(F) = -1$, then $u_1$ and $w_2$ are in one component, and $u_2$ and $w_1$ are in the other. This is identical to case 1, except the introduction of edges directed $u_1 \mapsto w_1$ and $u_2 \mapsto w_2$ form one cycle.

\textbf{Part II:} Let $c \in \widehat{\mathfrak{C}}^1_{\neq0}(L^0(G);u_1 u_2, w_1 w_2)$ and let $\check{c}\in {\mathfrak{C}}(L^0(G);u_1 u_2, w_1 w_2)$ be the unreduced contributor for $c$. Since $c$ is in a trivial activation class, $\check{c}$ must either (a) contain $2$ circles with $u_1 \mapsto w_1$ or $u_2 \mapsto w_2$ belonging to different circles, or (b) contain $1$ circle with $u_1 \mapsto w_1$ and $u_2 \mapsto w_2$ belonging to the same circle. 

\textit{Case 1 (Two-circles):} Suppose $\check{c}$ has exactly $2$-circles. First, unpack all backsteps of $c$, then re-introduce $u_1 \mapsto w_1$ and $u_2 \mapsto w_2$ to complete the two circles. Reverse the circle orientations and remove the adjacencies. The result is a Tutte-$2$-arborescence $F$ of the form $\langle u_1 w_1, u_2 w_2  \rangle$ and $sgn_T(F) = +1$. 

\textit{Case 2 (One-circle):} Again, this is similar to case 1, except the adjacencies introduced form a single circle. The result is a Tutte-$2$-arborescence $F$ of the form $\langle u_1 w_2, u_2 w_1  \rangle$, and $sgn_T(F) = -1$. \qed
\end{proof}

\begin{example}
To see how a Tutte-$2$-arborescence transforms into a circle-free reduced contributor, consider the top left Tutte-$2$-arborescence from Figure \ref{fig:2-arbors} in the calculation for $[v_5 v_4, v_6 v_1]$. This Tutte-$2$-arborescence appears on the left of Figure \ref{fig:TutteVsCont}. 

\begin{figure}[H]
    \centering
    \includegraphics{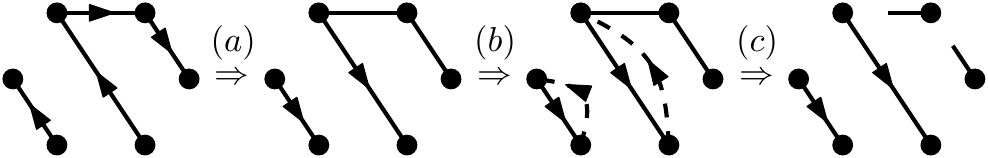}
    \caption{A Tutte-$2$-arborescence transforming into a reduced contributor.}
    \label{fig:TutteVsCont}
\end{figure}

The paths within each part of the arborescence are reversed in step $(a)$. The missing edge is added to produce a unique (directed) circle in step $(b)$. Next, all edges connected to each circle via a path are packed into backsteps away from each circle, producing the original restricted contributor. Finally, the introduced edges are removed to produce the reduced contributor in step $(c)$. 
\end{example}
We have the following immediate corollaries.

\begin{corollary}
Let $F$ be a Tutte-$2$-arborescence in the calculation of $[u_1 u_2, w_1 w_2]$ and $c_F$ be its corresponding element in $\widehat{\mathfrak{C}}^1_{\neq0}(L^0(G);u_1 u_2, w_1 w_2)$, then
\begin{enumerate}
    \item $sgn_T(F) = +1$ if, and only if, $\check{c}_F$ has exactly two cycles,
    \item $sgn_T(F) = -1$ if, and only if, $\check{c}_F$ has exactly one cycle.
\end{enumerate}
\end{corollary}

\begin{corollary} 
\label{c:unpackGivesTree}
Let $e$ be the edge between $w_1$ and $w_2$. Introducing the $w_1 w_2$-edge to any Tutte-$2$-arborescence associated to $[u_1 u_2, w_1 w_2]$ or a contributor-$2$-arborescence associated to an element of $\widehat{\mathfrak{C}}^1_{\neq0}(L^0(G);u_1 u_2, w_1 w_2)$ produces a spanning tree in $G \cup e$.
\end{corollary}

\begin{proof}
In either type of $2$-arborescence $w_1$ and $w_2$ are in different components and each component is a tree. If $e$ is an edge of $G$ a spanning tree of $G$ is produced. If $e$ does not exist in $G$, a spanning tree in $G \cup e$ is produced.
\qed \end{proof}

\begin{example}
Consider the graph in Figure \ref{fig:2-arbors}. Two of the reduced contributors in $\widehat{\mathfrak{C}}^1_{\neq 0}(L^0(G);v_5 v_4, v_6 v_1)$ that correspond to $[v_5 v_4, v_6 v_1]$ appear on the left of Figure \ref{fig:arbor-tree}. The middle figures are obtained by unpacking backsteps to produce a contributor-$2$-arborescence. Finally, the introduction of the $v_6 v_1$-edge yields a spanning tree.

\begin{figure}[H]
    \centering
    \includegraphics{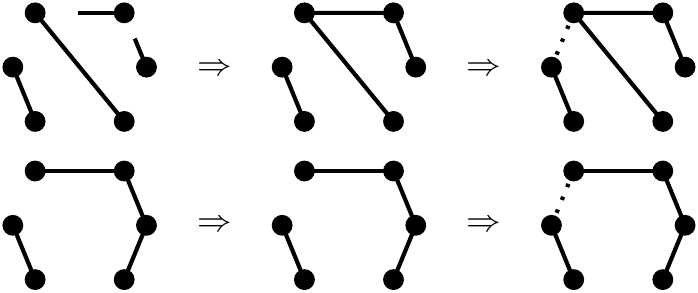}
    \caption{Trivial activation classes unpack into $2$-arborescences, and those used for edge labeling produce spanning trees.}
    \label{fig:arbor-tree}
\end{figure}

\end{example}

In the next Section Tutte-$2$-arborescences are replaced with entire contributor sets and the signs of the contributors and the non-trivial classes are characterized.

\section{Signed Graph Transpedances}
\label{sec:Transp}

We show that edge-labeling via signed contributors provide a generalization of transpedances and Kirchhoff-type Laws to signed graphs via the coefficients of the degree-$2$ monomials $x_{u_1w_1}x_{u_2w_2}$ from Theorem \ref{t:Poly}.

\subsection{Contributors as Transpedances}
\label{ssec:ContAreTransp}

The \emph{determinant-sign} of a contributor $c$ is taken from Theorem \ref{t:Poly}, where
\begin{align*}
    sgn_D(c) = (-1)^{ec(\check{c})+nc(c)+bs(c)}.
\end{align*}
The contributor-based transpedance for the determinant, or \emph{D-contributor-transpedance}, is defined as
\begin{align*}[u_1 u_2, w_1 w_2]_D = \dsum\limits_{c \in \widehat{\mathfrak{C}}_{\neq 0}(L^{0}(G);u_1 u_2, w_1 w_2) }
sgn_D(c),
\end{align*}
and consider the labeling of each $w_1 w_2$-edge with the signed contributors from $[u_1 u_2, w_1 w_2]_D$ when $w_1$ and $w_2$ are adjacent.

\begin{example}
Again, consider the graph in Figure \ref{fig:2-arbors}. The set of contributors that determine $[v_5 v_4, v_6 v_1]_D$, grouped into their activation classes, are shown in Figure \ref{fig:arbor-tree2}.

\begin{figure}[H]
    \centering
    \includegraphics{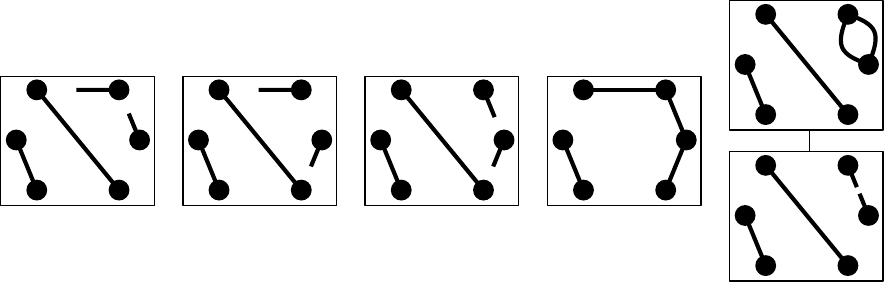}
    \caption{All activation classes for $[v_5 v_4, v_6 v_1]_D$.}
    \label{fig:arbor-tree2}
\end{figure}

We will see shortly that non-trivial classes sum to zero if all edges are positive, and only the trivial classes will determine $[v_5 v_4, v_6 v_1]_D$ if $G$ is a graph.
\end{example}

There is a simple relationship between the signs of a Tutte-$2$-arborescence and their associated reduced contributor.

\begin{lemma}
\label{l:TranspSignLemma}
Let $F$ be a Tutte-$2$-arborescence and $c_F$ be its corresponding element in $\widehat{\mathfrak{C}}^1_{\neq0}(L^0(G);u_1 u_2, w_1 w_2)$, then $sgn_T(F) = (-1)^{\lvert V \rvert} sgn_D(c_F)$.
\end{lemma}

\begin{proof} 

Tutte's transpedances $[u_1 u_2, w_1 w_2]$ are ordered second cofactors from the Laplacian $\mathbf{L}_{(G;u_1 u_2,w_1 w_2)}$, and the Tutte-$2$-arborescences are the signed commensurable parts that sum to $[u_1 u_2, w_1 w_2]$. From Theorem \ref{t:Poly}, the coefficient of $x_{u_1 w_1} x_{u_2 w_2}$ is $[u_1 u_2, w_1 w_2]_D$, and the reduced contributors are the signed commensurable parts that sum to $[u_1 u_2, w_1 w_2]_D$, but the coefficient of $x_{u_1 w_1} x_{u_2 w_2}$ is determined from $\mathbf{X} - \mathbf{L}_{G}$. The two adjacencies removed in each reduced contributor are mapped to $x_{u_1 w_1}$ and $x_{u_2 w_2}$, while all $\lvert V \rvert -2$ remaining Laplacian entries are negated; thus the sign discrepancy is $(-1)^{\lvert V \rvert -2} = (-1)^{\lvert V \rvert}$.
\qed \end{proof}

Tutte's transpedance degeneracy rule from Theorem \ref{t:TutteTransp} holds for $D$-contributor-transpedances. 

\begin{lemma}[Contributor Degeneracy]
\label{l:DDegeneracy}
Let $G$ be a signed graph with source $u_1$, sink $u_2$, and vertices $w_1$ and $w_2$, then 
\begin{align*}
    [u_1 u_1, w_1 w_2]_D = [u_1 u_2, w_1 w_1]_D = 0.
\end{align*}

\end{lemma}

\begin{proof}
The set of reduced contributors for $\widehat{\mathfrak{C}}_{\neq 0}(L^{0}(G);u_1 u_1, w_1 w_2)$ and $\widehat{\mathfrak{C}}_{\neq 0}(L^{0}(G);u_1 u_2, w_1 w_1)$ are both empty. The first would require two maps of the form $u_1 \mapsto w_1$ and $u_1 \mapsto w_2$, and there cannot be two tails at $u_1$. The second would require two maps of the form $u_1 \mapsto w_1$ and $u_2 \mapsto w_1$, and there cannot be two heads at $w_1$. \qed
\end{proof}

Since unreduced contributors represent permutation clones in $G$, we may apply the Linking Lemma on reduced contributors to produce $WU$-paths as the circles are cut. This shows that Tutte's energy reversal rule from Theorem \ref{t:TutteTransp} holds for $D$-contributor-transpedances.

\begin{lemma}[Contributor Energy Reversal]
\label{l:DReversal}
Let $G$ be a signed graph with source $u_1$, sink $u_2$, and vertices $w_1$ and $w_2$, then 
\begin{align*}
    [u_1 u_2, w_1 w_2]_D = -[u_1 u_2, w_2 w_1]_D = -[u_2 u_1, w_1 w_2]_D.
\end{align*}
\end{lemma}

\begin{proof}
We show the first equality, the second is similar.

Consider $c \in \widehat{\mathfrak{C}}_{\neq 0}(L^{0}(G);u_1 u_2, w_1 w_2)$ and reintroduce maps $u_1 \mapsto w_1$ and $u_2 \mapsto w_2$ to form the unreduced contributor $\check{c}$. There are two cases depending if $u_1, u_2, w_1, w_2$ belong to one or two circles in $\check{c}$.

\textit{Case 1 (Two circles):}
In this case we have that $\{u_1, w_1\}$ and $\{u_2, w_2\}$ are in disjoint circles in $\check{c}$. Remove $u_1 \mapsto w_1$ and $u_2 \mapsto w_2$ in $\check{c}$, and replace them with $u_1 \mapsto w_2$ and $u_2 \mapsto w_2$ to form a new non-zero unreduced contributor $\check{c}'$ in $\widehat{\mathfrak{C}}_{\neq 0}(L^{0}(G);u_1 u_2, w_2 w_1)$ where $\{u_1,u_2, w_1, w_2\}$ are in a single circle. Since $c$ and $c'$ have the same adjacencies and backsteps, the sign difference between $sgn_D(c)$ and $sgn_D(c')$ is determined by the even circle structure of their unreduced contributors. If the original circles were both even, the new single circle is even; a loss of one even circle. If the original circles were both odd, the new circle is even; a gain of one even circle. If the original circles have different parity, the new circle is odd; a loss of one even circle. In any case $sgn_D(c) = -sgn_D(c')$.

\textit{Case 2 (One circle):}
In this case we have that $\{u_1,u_2, w_1, w_2\}$ are in a single circle $\check{c}$. Remove $u_1 \mapsto w_1$ and $u_2 \mapsto w_2$ in $\check{c}$, and replace them with $u_1 \mapsto w_2$ and $u_2 \mapsto w_2$ to form a new non-zero unreduced contributor $\check{c}'$ in $\widehat{\mathfrak{C}}_{\neq 0}(L^{0}(G);u_1 u_2, w_2 w_1)$ where $\{u_1, w_1\}$ and $\{u_2, w_2\}$ are in disjoint circles. Since $c$ and $c'$ have the same adjacencies and backsteps, the sign difference between $sgn_D(c)$ and $sgn_D(c')$ is determined by the even circle structure of their unreduced contributors. If the original circle is even, each new circle is odd. If the original circle is odd, each new circle is even.

This process is reversible, so we have the first equality. The second equality is similar. \qed
\end{proof}

\subsection{Transpedance Evaluation}
\label{ssec:TranpEval}

With the sign adjustment in Lemma \ref{l:TranspSignLemma}, we can immediately use the elements of $\widehat{\mathfrak{C}}^1_{\neq0}(L^0(G);u_1 u_2, w_1 w_2)$ in place of transpedances. However, we now characterize the placement of entire contributor families on edges. Let $\hat{\mathfrak{A}}(\mathbf{u},\mathbf{w};G)$ be the set of all reduced activation classes of the form $\hat{\mathcal{A}}(\mathbf{u},\mathbf{w};G)$, and let $\hat{\mathfrak{A}}^-(\mathbf{u},\mathbf{w};G)$ be the subset of $\hat{\mathfrak{A}}(\mathbf{u},\mathbf{w};G)$ such that no element contains a positive circle. The activation class transversal consisting of maximal elements is denoted by $\mathcal{M}_{\mathbf{u},\mathbf{w}}$, and $\mathcal{M}_{\mathbf{u},\mathbf{w}}^-$ is the subset of maximal elements that are positive-circle-free.
Since each activation class is Boolean (Lemma \ref{l:booleanLattice}), $D$-contributor-transpedances have a simple presentation via the maximal element of each activation class. 

\begin{theorem}
\label{t:signingBT} 
If $G$ is a signed graph, then 
\begin{align*}
    [u_1u_2,w_1w_2]_D = \dsum\limits_{m \in \mathcal{M}_{(u_1u_2,w_1w_2)}^-} sgn_D(m) \cdot (2)^{\eta{(m)}}
\end{align*}
where $\eta{(m)}$ is the number of negative circles in maximal contributor $m$.
\end{theorem}

\begin{proof}

Let $G$ be a signed graph with distinguished source $u_1$, sink $u_2$, edge $w_1 w_2$, and total orderings $\mathbf{u}=(u_1,u_2)$ and $\mathbf{w}=(w_1,w_2)$. Also, let $\hat{\mathfrak{A}} = \hat{\mathfrak{A}}(\mathbf{u},\mathbf{w};G)$, and $\hat{\mathcal{A}} = \hat{\mathcal{A}}(\mathbf{u},\mathbf{w};G)$.
Partition the $D$-contributor-transpedance value $[u_1 u_2, w_1 w_2]_D$ into activation classes as follows:
\begin{align*}[u_1 u_2, w_1 w_2]_D  &=\dsum\limits_{ \hat{\mathcal{A}} \in \hat{\mathfrak{A}}}\dsum\limits_{c \in \hat{\mathcal{A}}} sgn_D(c) \\ 
&= \dsum\limits_{\hat{\mathcal{A}} \in \hat{\mathfrak{A}}^-}\dsum\limits_{c \in \hat{\mathcal{A}} } sgn_D(c) + \dsum\limits_{ \hat{\mathcal{A}} \in \hat{\mathfrak{A}} \setminus \hat{{\mathfrak{A}}}^- }\dsum\limits_{c \in \hat{\mathcal{A}}} sgn_D(c)
\end{align*}

From Lemma \ref{l:booleanLattice}, activation classes form Boolean lattices, and each sum is calculated separately. 

\textit{Case 1 (No positive circles):} Let contributors $c$ and $c'$ only differ by a single negative circle, which appears in $c$ but not in $c'$. Let the length of this circle be $\ell$. Packing this circle into backsteps will yield a loss of a single positive circle and a gain of $\ell$ backsteps.

\textit{Case 1a ($\ell$ is odd):} If $\ell$ is odd, the $sgn_D(c')$ is related to $sgn_D(c)$ as follows:
\begin{align*}
     sgn_D(c') &= (-1)^{ec(\check{c})+(nc(c)-1)+(bs(c)+\ell)} \\
    &= (-1)^{ec(\check{c})+nc(c)+bs(c)}\cdot (-1)^{\ell-1} = sgn_D(c).
\end{align*}

\textit{Case 1b ($\ell$ is even):} If $\ell$ is even, packing also loses an even circle and $sgn_D(c')$ is related to $sgn_D(c)$ as follows:
\begin{align*}
    sgn_D(c') &= (-1)^{(ec(\check{c})-1)+(nc(c)-1)+(bs(c)+\ell)} \\
    &= (-1)^{ec(\check{c})+nc(c)+bs(c)}\cdot (-1)^{\ell-2} = sgn_D(c).
\end{align*}

Since each element has the same sign and each activation class is Boolean there are $2^{\,\eta{(m)}}$ contributors, where $m$ is the maximal contributor containing $\eta{(m)}$ circles.

\textit{Case 2 (Positive circle):} Let  contributors $c$ and $c'$ only differ by a single positive circle, which appears in $c$ but not in $c'$. Let the length of this circle be $\ell$. Packing this circle into backsteps will yield a loss of a single positive circle and a gain of $\ell$ backsteps.

\textit{Case 2a ($\ell$ is odd):} If $\ell$ is odd, the $sgn_D(c')$ is related to $sgn_D(c)$ as follows:
\begin{align*}
     sgn_D(c') &= (-1)^{ec(\check{c})+nc(c)+(bs(c)+\ell)} \\
    &= (-1)^{ec(\check{c})+nc(c)+bs(c)}\cdot (-1)^{\ell} = -sgn_D(c).
\end{align*}

\textit{Case 2b ($\ell$ is even):}
If $\ell$ is even, packing also loses an even circle and $sgn_D(c')$ is related to $sgn_D(c)$ as follows:
\begin{align*}
     sgn_D(c') &= (-1)^{(ec(\check{c})-1)+nc(c)+(bs(c)+\ell)} \\
    &= (-1)^{ec(\check{c})+nc(c)+bs(c)}\cdot (-1)^{\ell - 1} = -sgn_D(c).
\end{align*}

Again, since each activation class is Boolean, there is a bijection between contributors with circle $C$ and those without $C$ via packing/unpacking. Thus, each activation class that contains a contributor with a positive circle will have those contributors sum to $0$. Moreover, the remaining classes are determined by the sign of their maximal element.
\begin{align*}[u_1 u_2, w_1 w_2]_D  
&= \dsum\limits_{\hat{\mathcal{A}} \in \hat{\mathfrak{A}}^-}\dsum\limits_{c \in \hat{\mathcal{A}} } sgn_D(c) + 0 \\ 
&= \dsum\limits_{m \in \mathcal{M}_{\mathbf{u};\mathbf{w}}^-} sgn_D(m) \cdot (2)^{\eta{(m)}}.
\end{align*}
\qed \end{proof}

Combining Lemmas \ref{l:ArborBijection}, \ref{l:TranspSignLemma}, and Theorem \ref{t:signingBT} we have the following interpretation of Tutte-transpedances:

\begin{corollary}[Parity-Polarity Reversal]
\label{c:IsTutte}
If $G$ is a signed graph with all positive edges, then $[u_1 u_2, w_1 w_2] = (-1)^{\lvert V \rvert}[u_1 u_2, w_1 w_2]_D$. 
\end{corollary}

\begin{proof}
If all edges are positive, by Theorem  \ref{t:signingBT}, the only non-cancellative terms are trivial reduced activation classes. The bijection between $2$-arborescence types in Lemma \ref{l:ArborBijection} combined with the signing in Lemma \ref{l:TranspSignLemma} completes the proof. \qed
\end{proof}

That is, Tutte's edge-labeling via transpedances provides a natural orientation from source to sink, while the contributor version is reversed for graphs with an odd number of vertices. 

\begin{example}
If all edges are positive, the contributors for $[v_5 v_4, v_6 v_1]_D$ in Figure \ref{fig:arbor-tree2} produce a value of $+4$ as the non-trivial classes sum to $0$ and there is an even number of vertices. This agrees with Tutte's $[v_5 v_4, v_6 v_1]$. 
\end{example}

However, a signed graph may not have their non-trivial activation classes cancel.

\begin{example} 
For a new example, consider the signed graph in Figure \ref{fig:NegativeCycleT} with source $v_5$ and sink $v_4$. To calculate the D-contributor-transpedance along edge $v_5 v_4$ we examine $[v_5 v_4, v_5 v_4]_D$. The contributors in Figure \ref{fig:NegativeCycleT} are the non-cancellative contributors as they do not contain positive circles. Since there are an odd number of vertices the value is negated relative to Tutte's and the value is $-12$.
\begin{figure}[H]
    \centering
    \includegraphics{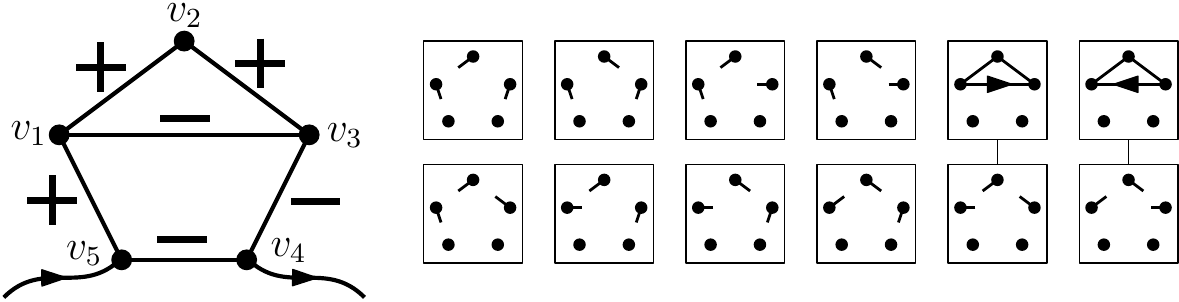}
    \caption{Non-trivial reduced-contributors signed $[v_5 v_4, v_5 v_4]_D = -12$}
    \label{fig:NegativeCycleT}
\end{figure}
If all the edges were positive, the $D$-contributor-transpedance would have been $-8$ as the oriented $3$-circles would be positive. Also, observe that there are far more than $12$ contributors on edge $v_5 v_4$, as the contributors in Figure \ref{fig:NegativeCycleT2} always cancel as they repeat an adjacency, so the circle is always positive.
\begin{figure}[H]
    \centering
    \includegraphics{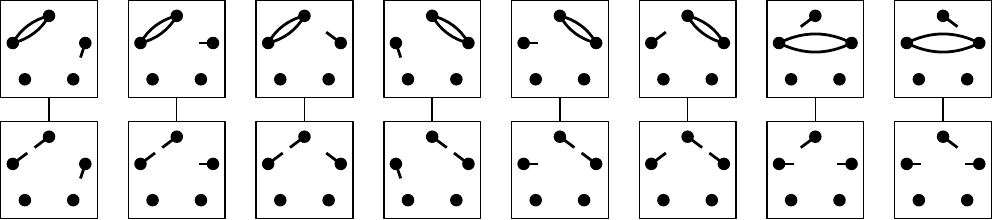}
    \caption{Non-trivial reduced-contributors signed $[v_5 v_4, v_5 v_4]_D = -12$}
    \label{fig:NegativeCycleT2}
\end{figure}

The sign between an edge does not matter when determining the $D$-contributor-transpedance on that edge, as the edge fails to appear in any contributor. Additionally, cycles may not cancel in their activation class as they would in a graph.
\end{example}

\section{Adjacency Exchange and Kirchhoff's Laws}
\label{sec:Kirch}

Contributor-transpedances satisfy their own general Degeneracy (Lemma \ref{l:DDegeneracy}) and Energy Reversal (Lemma \ref{l:DReversal}) Kirchhoff-type laws. They are evaluated via activation classes (Theorem \ref{t:signingBT}), with all positive graphs related to Tutte-transpedances via Polarity reversal (Corollary \ref{c:IsTutte}). The Cycle Conservation and Vertex Conservation properties from Theorem \ref{t:TutteTransp} are now investigated by showing that transpedances are contributor sorting along source-sink paths as a generalization of Corollary \ref{c:unpackGivesTree}. However, the expectation of conservation cannot be expected.

\subsection{Source-sink Pathing}

Let $G$ be a signed graph with source $u_1$, sink $u_2$, and vertices $w_1$ and $w_2$. If $w_1$ and $w_2$ are adjacent, call their edge $e$. If $w_1$ and $w_2$ are not adjacent, regard $G$ as a subgraph $G \cup e_{w_1 w_2}$ where edge $e_{w_1 w_2}$ is added between $w_1$ and $w_2$. This is called the \emph{local-loading of $G$ at $\{w_1, w_2\}$}, and is related to the injective loading properties from \cite{Grill1,IH2}; to simplify notation we will simply write $G \cup e_{w_1 w_2}$ with the understanding that $e_{w_1 w_2}$ may exist in $G$. Let $\mathcal{P}(u_1 u_2, w_1 w_2)$ be the set of $u_1 u_2$-paths containing $e_{w_1 w_2}$ in $G \cup e_{w_1 w_2}$.

\begin{lemma}
\label{l:PathLem}
A contributor $c$ is in $\widehat{\mathfrak{C}}_{\neq 0}(L^{0}(G);u_1 u_2, w_1 w_2)$ if, and only if, $c$ contains a unique path $P \in \mathcal{P}(u_1 u_2, w_1 w_2)$ in $G \cup e_{w_1 w_2}$. 
\end{lemma}

\begin{proof}

\textbf{Part I:} Let $c \in \widehat{\mathfrak{C}}_{\neq 0}(L^{0}(G);u_1 u_2, w_1 w_2)$, and $\check{c}$ be its unreduced contributor. There are two cases depending if $u_1, u_2, w_1, w_2$ belong to one or two circles in $\check{c}$.

\textit{Case 1 (Two circles):}
In this case we have $u_1, w_1$ and $u_2, w_2$ are in disjoint circles in $\check{c}$. Remove $u_1 \mapsto w_1$ and $u_2 \mapsto w_2$, reverse $u_1 \rightarrow w_1$, and introduce edge $e$ to produce a $u_1 u_2$-path $P$ where $w_1$ precedes $w_2$ in $P$. Any additional circles and backsteps are external and may only extend the activation class. 

\textit{Case 2 (One circle):}
In this case we have $u_1, u_2, w_1, w_2$ in a single circle in $\check{c}$. Remove $u_1 \mapsto w_1$ and $u_2 \mapsto w_2$, reverse $u_1 \rightarrow w_2$, and introduce edge $e$ to produce a $u_1 u_2$-path $P$ where $w_2$ precedes $w_1$ in $P$. Any additional circles and backsteps are external and may only extend the activation class. 

\textbf{Part II:} Let $P \in \mathcal{P}(u_1 u_2, w_1 w_2)$. There are two cases depending if $w_1$ precedes $w_2$ or $w_2$ precedes $w_1$ in $P$.

\textit{Case 1 ($w_1$ precedes $w_2$):}
Delete $e_{w_1 w_2}$, and introduce $u_1 \mapsto w_1$ and $u_2 \mapsto w_2$. Reverse the $u_1 w_1$-part of $P$, and do not reverse the $u_2 w_2$-part of $P$ to make two circles. Introduce backsteps/circles at all remaining vertices for form an unreduced contributor $\check{c}$. Remove $u_1 \mapsto w_1$ and $u_2 \mapsto w_2$ to get $c \in \widehat{\mathfrak{C}}_{\neq 0}(L^{0}(G);u_1 u_2, w_1 w_2)$. Pack/unpack as necessary to form activation classes.  

\textit{Case 2 ($w_2$ precedes $w_1$):}
Delete $e_{w_1 w_2}$, and introduce $u_1 \mapsto w_1$ and $u_2 \mapsto w_2$. Reverse the $u_1 w_2$-part of $P$, and do not reverse $u_2 w_1$-part of $P$ to make one circle. Introduce backsteps/circles at all remaining vertices to form an unreduced contributor $\check{c}$. Remove $u_1 \mapsto w_1$ and $u_2 \mapsto w_2$ to get $c \in \widehat{\mathfrak{C}}_{\neq 0}(L^{0}(G);u_1 u_2, w_1 w_2)$. Pack/unpack as necessary to form activation classes. \qed
\end{proof}

We now have the immediate corollaries demonstrating that all contributors for a given transpedance are related to direct source-sink path property. 

\begin{corollary}
\label{c:ShortestPath}
Let $c \in \widehat{\mathfrak{C}}_{\neq 0}(L^{0}(G);u_1 u_2, w_1 w_2)$. Every edge-adjacency appearing in $c$ outside an activated circle is in one of the parts of the $w_i u_j$-paths. Moreover, these paths are oriented from $w_i$ to $u_j$.
\end{corollary}

\begin{corollary}
\label{c:PendantZero}
If $w_1$ is a monovalaent vertex that is not a source or sink with supporting edge $e_{w_1 w_2}$, then $[u_1 u_2, w_1 w_2] = 0$.
\end{corollary}

We also have the following simple interpretation of Tutte's transpedances.

\begin{corollary}
\label{c:TreeSort}
Let $G$ be a graph with source $u_1$ and sink $u_2$. The edge labeling of $G$ by transpedances $[u_1 u_2, w_1 w_2]$ is equivalent to a sorting of spanning trees via adjacency swapping along the $u_1 u_2$-path in $G \cup e$.
\end{corollary}

\begin{proof}
Corollary \ref{c:unpackGivesTree} provides an interpretation of trivial activation classes as spanning trees, even for transpedances not on adjacencies. Additionally, Corollary \ref{c:IsTutte} shows that these are the only objects that survive cancellation in the Boolean activation classes in a graph. Part 4 of Theorem \ref{t:TutteTransp} indicates the net inflow and outflow is the tree-number. \qed
\end{proof}

\begin{example}
The reduced contributors in trivial activation classes for Figure \ref{fig:KirchA} appear on each edge in Figure \ref{fig:ExampleA} (left). A source-sink path is indicated on the right, and the associated unpacked contributors appear with edge inserted to produce spanning trees to visualize both the spanning tree sorting and the unique pathing property.

\begin{figure}[H]
    \centering
    \includegraphics{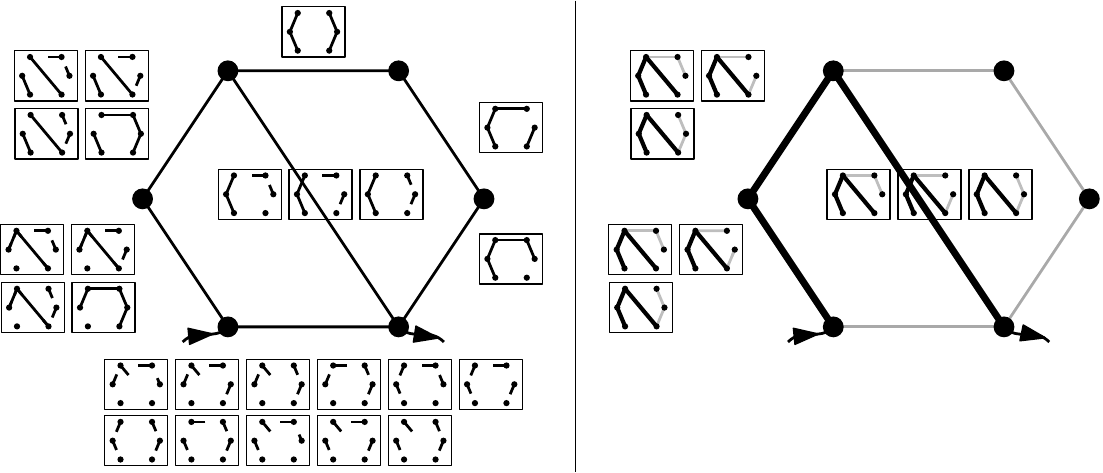}
    \caption{Left: Contributors from trivial classes. Right: The associated spanning trees and unique paths.}
    \label{fig:ExampleA}
\end{figure}

\end{example}

Combining Theorem \ref{t:signingBT} and Lemma \ref{l:PathLem} we can discuss the Cycle and Vertex Conservation properties from Theorem \ref{t:TutteTransp}. These conservation laws are byproducts of non-trivial Boolean classes vanishing in a graph, coupled with a natural matching of the elements in the trivial classes whose signs cancel. Negative edges may produce non-vanishing Boolean classes as well as matched trivial classes of the same sign.

\begin{lemma}[Contributor Cycle ``Conservation'']
\label{l:DCycle}
There is a matching between the elements of 
\begin{align*}
\widehat{\mathfrak{C}}^1_{\neq 0}(L^{0}(G);u_1 u_2, w_1 w_2) \cup \widehat{\mathfrak{C}}^1_{\neq 0}(L^{0}(G);u_1 u_2, w_2 w_3) \cup \widehat{\mathfrak{C}}^1_{\neq 0}(L^{0}(G);u_1 u_2, w_3 w_1).
\end{align*}
\end{lemma}

\begin{proof}

Let $G$ be a signed graph with source $u_1$, sink $u_2$, and vertices $w_1$, $w_2$, and $w_3$. Additionally, let $e_{w_1 w_2}$, $e_{w_2 w_3}$, and $e_{w_3 w_1}$ be the edges between their respective vertices, or the edge introduced to $G$ if one does not exist. 

Consider $c \in \widehat{\mathfrak{C}}^1_{\neq 0}(L^{0}(G);u_1 u_2, w_1 w_2)$. From Lemma \ref{l:PathLem}, let $P$ be the unique $u_1 u_2$-path in $c$ made with the inclusion of $e_{w_1 w_2}$ so that $P = P_{u_1, w_i} \cup e_{w_1 w_2} \cup P_{w_j, u_2}$, where $\{i,j\} = \{1,2\}$. Since $c$ is in a trivial class, there are no circles to activate, and from Corollary \ref{c:ShortestPath} vertex $w_3$ must be linked to $P_{u_1, w_i}$ or $P_{w_j, u_2}$ by a sequence of unpackings. Moreover, all backsteps outside of circle-activation unpack towards $P$, so there is a unique vertex $w'$ that meets exactly one of  $P_{u_1, w_i}$ or $P_{w_j, u_2}$.

\begin{center}
\includegraphics{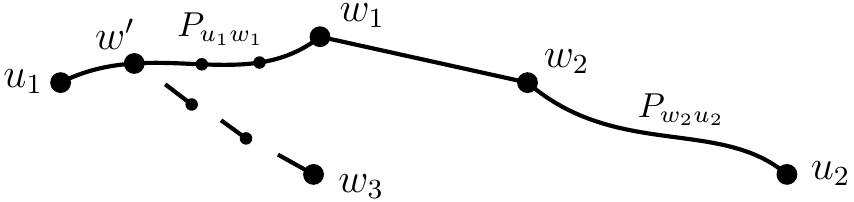}
\end{center}

Assume $w'$ meets $P_{u_1, w_i}$, the case where $w'$ meets $P_{w_j, u_2}$ is similar. Form the path $P': u_1 \rightarrow w' \rightarrow w_3$. $P'$ may contain $w_i$ if $w' = w_i$ but cannot contain $e_{w_1 w_2}$. Introducing edge $e_{w_3,w_j}$ forms a unique $u_1 u_2$-path $P'' = P' \cup e_{w_3,w_j} \cup P_{w_j, u_2}$ that uses exactly one of $e_{w_1 w_2}$, $e_{w_2 w_3}$, or $e_{w_3 w_1}$. Removing $e_{w_3,w_i}$, reversing the $u_1 w_i$ part of $P''$, and packing all non-$P''$ adjacencies away from $P''$ leaves a unique contributor $c' \in \widehat{\mathfrak{C}}^1_{\neq 0}(L^{0}(G);u_1 u_2, w_3 w_j)$. 

\begin{center}
\includegraphics{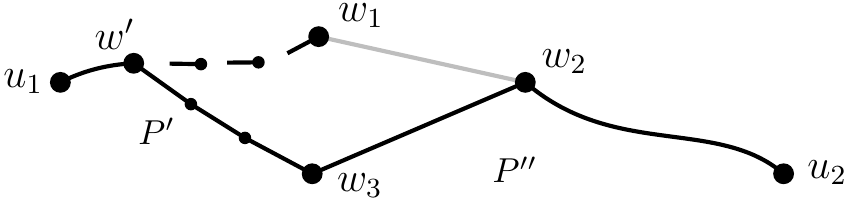}
\end{center}

Moreover, there is no corresponding contributor in $\widehat{\mathfrak{C}}^1_{\neq 0}(L^{0}(G);u_1 u_2, w_3 w_i)$ as $e_{w_3,w_i}$ does not form a path without using more than one of $e_{w_1 w_2}$, $e_{w_2 w_3}$, and $e_{w_3 w_1}$. \qed

\end{proof}

Tutte's Cycle Conservation in Theorem \ref{t:TutteTransp} is an immediate consequence of Lemma \ref{l:DCycle} as every edge is positive and the matching converts between one and two circles, changing their signs. General signed graphic conservation, however, cannot be guaranteed as (1) there may be negative edges between a trivial-class matching, and (2) the non-trivial classes need not cancel. Tutte's Vertex Conservation in Theorem \ref{t:TutteTransp} is also an immediate consequence of the following lemma and is easily seen in Figure \ref{fig:ExampleA} by following the contributor sorting along source-sink paths. 

\begin{lemma}[Vertex ``Conservation'']
\label{l:DVertex}
Let $G$ be a signed graph with source $u_1$, sink $u_2$, and let $v$ be another vertex. 
\begin{enumerate}
    \item If $v \nin \{u_1 , u_2\}$, then $\displaystyle\lvert \bigcup_{x \sim v}\widehat{\mathfrak{C}}^1_{\neq 0}(L^{0}(G);u_1 u_2, x v) \rvert = \lvert \bigcup_{y \sim v}\widehat{\mathfrak{C}}^1_{\neq 0}(L^{0}(G);u_1 u_2, v y) \rvert$.
    \item If $v \in \{u_1 , u_2\}$, then $\displaystyle\lvert \bigcup_{x \sim v}\widehat{\mathfrak{C}}^1_{\neq 0}(L^{0}(G);u_1 u_2, u_1 x) \rvert = \lvert \bigcup_{y \sim v}\widehat{\mathfrak{C}}^1_{\neq 0}(L^{0}(G);u_1 u_2, y u_2) \rvert$.
\end{enumerate}

\end{lemma}
\begin{proof}

Let $G$ be a signed graph with source $u_1$, sink $u_2$, and let $v$ be another vertex. Consider $\widehat{\mathfrak{C}}^1_{\neq 0}(L^{0}(G);u_1 u_2, x v)$ and $\widehat{\mathfrak{C}}^1_{\neq 0}(L^{0}(G);u_1 u_2, v y)$, where the edges $e_{x v}$ and $e_{v y}$ exist in $G$.

\textit{Case 1 ($v \nin \{u_1 , u_2\}$):} If $c \in \widehat{\mathfrak{C}}^1_{\neq 0}(L^{0}(G);u_1 u_2, x v)$, using Lemma \ref{l:PathLem} consider the $u_1 u_2$-paths that contains $e_{x v}$. Since $v$ is not the source or sink, each path must contain exactly one of the edges $e_{v y}$ for some $y$. From Corollary \ref{c:ShortestPath}, all contributors in $\widehat{\mathfrak{C}}^1_{\neq 0}(L^{0}(G);u_1 u_2, x v)$ associated to a path containing both $e_{x v}$ and $e_{v y}$ must also have a corresponding element in $\widehat{\mathfrak{C}}^1(L^{0}(G);u_1 u_2, v y)$.

The argument is identical on the preceding edge when starting with $\widehat{\mathfrak{C}}^1(L^{0}(G);u_1 u_2, v y)$.

\textit{Case 2 ($v \in \{u_1 , u_2\}$):}
If $v$ is the source, there are no $v$-entrant edges in any $u_1 u_2$-path. While, if $v$ is the sink, there are no $v$-salient edges in any $u_1 u_2$-path. However, from Lemma \ref{l:PathLem} and Corollary \ref{c:ShortestPath}, all contributors arise from $u_1 u_2$-paths, therefore all trivial class contributors out of $u_1$ have a corresponding contributor in to $u_2$. \qed
\end{proof}

If every edge of a signed graph is positive, not only do the non-trivial activation classes sum to zero, but the trivial ones in each matching above also cancel. Thus, conservation is guaranteed when $G$ has all positive edges.

\begin{corollary}
If $G$ has all positive edges, then the D-contributor-transpedances are both cycle and vertex conservative.
\end{corollary}

It is clear that a graph with a single negative edge that is also between the source and sink is conservative, as that edge never appears in any contributor. It seems worthwhile to produce a complete characterization of signed graphs that are conservative, even for fixed source and sinks.

\section{Maximizing Transpedance, Permanents, and Signless Laplacians} 
\label{sec:MaxPerm}

\subsection{Oriented Hypergraphs and Contributor Counting}
\label{ssec:ContCount}

Since contributor mappings are used in determining the characteristic and total minor polynomials \cite{OHSachs,IH2}, and generalizations of the Matrix-Tree Theorem \cite{OHMTT}, we examine the net placement of contributors on a given graph. The permanent of the oriented hypergraphic signless Laplacian was shown to count the number of contributors, which occurs when every adjacency is negative.

\begin{theorem}[\cite{OHSachs}, Theorem 4.3.1 part 1]
Let $G$ be an oriented hypergraph with no isolated vertices or $0$-edges with Laplacian matrix $\mathbf{L}_{G}$, then $\mathrm{perm}(\mathbf{L}_{G})=\left\vert \mathfrak{C}(G)\right\vert $ if, and only if, every edge of $G$ is extroverted or introverted.
\label{OHSachs2}
\end{theorem}

As in prior sections, the previous theorem was a direct calculation on the Laplacian, while we make use of the coefficient of the total minor polynomial to keep track of the ordered minor placement. The \emph{permanental-sign} of a contributor $c$ is taken from Theorem \ref{t:Poly}, where \begin{align*}
    sgn_P(c) = (-1)^{nc(c)+bs(c)}.
\end{align*}
The signless Laplacian can be used to count the number of reduced contributors for any oriented hypergraph.

\begin{theorem}
\label{t:naInductiveT}
If $G$ is an oriented hypergraph with all negative adjacencies, then 
\begin{align*}
    \dsum\limits_{c \in \widehat{\mathfrak{C}}_{\neq 0}(L^{0}(G);\mathbf{u},\mathbf{w}) } sgn_P(c) = (-1)^{\lvert V \rvert - k} \; \lvert \widehat{\mathfrak{C}}_{\neq 0}(L^{0}(G);\mathbf{u},\mathbf{w}) \rvert
\end{align*}
\end{theorem}

\begin{proof}
Let $k = \vert U \rvert = \vert W \rvert$ with $U$ and $W$ totally orderings $\mathbf{u}$ and $\mathbf{w}$. Also let $sgn_P(c) = (-1)^{nc(c)\,+\,bs(c)}$ be the permanent signing function.
We proceed with an inductive argument: 

\textit{Case 1 ($k=0$):} Observe that if $c$ is a minimal (identity-clone) contributor, then $nc(c)=0$ and $bs(c) = \lvert V \rvert$, and the permanent sign of all minimal contributors is $(-1)^{\lvert V \rvert}$. 
If $c'$ is any contributor that can unpack into another covering contributor $c''$ containing a new cycle of length $\ell$, then we have two cases based on $\ell$'s parity.

\textit{Case 1a ($\ell$ is even):} Unpack $\ell$ backsteps in $c'$ to form a cycle of length $\ell$ in $c''$. Since $\ell$ is even and all edges are negative, we lose $\ell$ backsteps and gain $0$ negative components. Since $-\ell+0$ is even, $sgn_P(c') = sgn_P(c'')$.

\textit{Case 1b ($\ell$ is odd):}
Unpack $\ell$ backsteps in $c'$ to form a cycle of length $\ell$ in $c''$. Since $\ell$ is odd and all edges are negative, we lose $\ell$ backsteps and gain $1$ negative component. Since $-\ell+1$ is even, $sgn_P(c') = sgn_P(c'')$.

Thus, all contributors have the same sign as their minimal contributor, and all minimal contributors have the same sign, giving
\begin{align*}
    \dsum\limits_{c \in \widehat{\mathfrak{C}}_{\neq 0}(L^{0}(G);\mathbf{u},\mathbf{w}) } sgn_P(c) = (-1)^{\lvert V \rvert} \; \lvert\widehat{\mathfrak{C}}(G) \rvert.
\end{align*}

\textit{Case $2$ ($k>0$):}
In a contributor with all negative adjacencies, deleting a negative edge will swap the sign of the component that contained the edge, thus changing $nc(c)$ by one. Deleting a backstep will decrease $bs(c)$ by one, which will flip the permanent sign of the total contributor. Since all contributors in $\widehat{\mathfrak{C}}(G)$ have the same permanent signing,
the sign alternates with every edge or backstep that is removed. Thus, the permanent counts of reduced contributors must be $(-1)^{\lvert V \rvert-k} |\widehat{\mathfrak{C}}_{\neq 0}(L^{0}(G);\mathbf{u},\mathbf{w})|$.
\qed \end{proof}

\subsection{Signed Graphs and Maximal Transpedance}
\label{ssec:MaxTransp}

We define the contributor based transpedance for the permanent, or \emph{P-contributor-transpedance}, to be
\begin{align*}[u_1 u_2, w_1 w_2]_P = \dsum\limits_{c \in \widehat{\mathfrak{C}}_{\neq 0}(L^{0}(G);u_1 u_2, w_1 w_2) }
sgn_P(c).
\end{align*}

\begin{corollary}
If $G$ is an oriented hypergraph with all negative adjacencies, then 
\begin{align*}
    [u_1 u_2, w_1 w_2]_P = (-1)^{\lvert V \rvert} \; \lvert \widehat{\mathfrak{C}}_{\neq 0}(L^{0}(G);u_1 u_2, w_1 w_2) \rvert
\end{align*}
\end{corollary}

\begin{proof}
From Theorem \ref{t:naInductiveT} the value of $[u_1 u_2, w_1 w_2]_P$ is equal to $(-1)^{\lvert V \rvert - 2}$ $\lvert \widehat{\mathfrak{C}}_{\neq 0}(L^{0}(G);u_1 u_2, w_1 w_2) \rvert$, and $\lvert V \rvert - 2$ and $\lvert V \rvert$ have the same parity. \qed 
\end{proof}

While a contributor-type transpedance can be defined for an arbitrary oriented hypergraph, the parallel adjacencies and the non-Boolean structure for tail-equivalence classes require further examination. However, the signless Laplacian provides a count on the total number of contributors, and a maximum possible value for a bidirected graph.  

A permanent version of Tutte's transpedance Theorem follows. 

\begin{lemma}[P-Contributor Degeneracy]
\label{l:PKirch}
Let $G$ be a signed graph with source $u_1$, sink $u_2$, and distinct vertices $w_1$, $w_2$ and $w_3$, then 
\begin{enumerate}
    \item $[u_1 u_1, w_1 w_2]_P = [u_1 u_2, w_1 w_1]_P = 0$,
    \item $[u_1 u_2, w_1 w_2]_P = [u_1 u_2, w_2 w_1]_P = [u_2 u_1, w_1 w_2]_P$,
    \item There is a matching between the elements of 
        \begin{align*}
        \widehat{\mathfrak{C}}^1_{\neq 0}(L^{0}(G);u_1 u_2, w_1 w_2) \cup \widehat{\mathfrak{C}}^1_{\neq 0}(L^{0}(G);u_1 u_2, w_2 w_3) \cup \widehat{\mathfrak{C}}^1_{\neq 0}(L^{0}(G);u_1 u_2, w_3 w_1).
        \end{align*}
    \item Let $G$ be a signed graph with source $u_1$, sink $u_2$, and let $v$ be another vertex. 
\begin{enumerate}
    \item If $v \nin \{u_1 , u_2\}$, then $\displaystyle\lvert \bigcup_{x \sim v}\widehat{\mathfrak{C}}^1_{\neq 0}(G;u_1 u_2, x v) \rvert = \lvert \bigcup_{y \sim v}\widehat{\mathfrak{C}}^1_{\neq 0}(G;u_1 u_2, v y) \rvert$.
    \item If $v \in \{u_1 , u_2\}$, then $\displaystyle\lvert \bigcup_{x \sim v}\widehat{\mathfrak{C}}^1_{\neq 0}(G;u_1 u_2, u_1 x) \rvert = \lvert \bigcup_{y \sim v}\widehat{\mathfrak{C}}^1_{\neq 0}(G;u_1 u_2, y u_2) \rvert$.
\end{enumerate}
\end{enumerate}
\end{lemma}
\begin{proof}
The proofs are identical to the determinant case as they are the same set of objects. The only exception is even circles are not included in any signs. \qed
\end{proof}

\section*{Acknowledgments}
A portion of this research was supported by Mathworks at Texas State University. The authors would like to thank Max Warshauer and Amelia Hu for their support and discussions.

\bibliographystyle{plain}

\end{document}